\documentclass[12pt]{amsart}

\usepackage{amsfonts,amssymb,amsmath,graphicx,hyperref}

\numberwithin{equation}{section}
\newtheorem{theorem}{Theorem}[section]
\newtheorem{lemma}{Lemma}[section]
\newtheorem{prop}{Proposition}[section]

\begin{document}

\title{The AJ Conjecture for Cables of Two-Bridge Knots}

\author{Nathan Druivenga}
\address{Mathematics Department\\
University of Iowa\\
Iowa City, IA}
\email{nathan-druivenga@uiowa.edu}

\begin{abstract}
The $AJ$-conjecture for a knot $K \subset S^3$ relates the $A$-polynomial and the colored Jones polynomial of $K$.  If a two-bridge knot $K$ satisfies the $AJ$-conjecture, we give sufficient conditions on $K$ for the $(r,2)$-cable knot $C$ to also satisfy the $AJ$-conjecture.  If a reduced alternating diagram of $K$ has $\eta_+$ positive crossings and $\eta_-$ negative crossings, then $C$ will satisfy the $AJ$-conjecture when $(r+4\eta_-)(r-4\eta_+)>0$ and the conditions of the main theorem are satisfied. 
\end{abstract}
\maketitle
\vspace{-.9 cm}
\section{Introduction}
\label{sec:intro}

The $AJ$-conjecture \cite{Ga} is a proposed relationship between two different invariants of a knot $K$ in $S^3$, the $A$-polynomial and the colored Jones polynomial.  The $A$-polynomial is determined by the fundamental group of the knot complement, while the colored Jones polynomial is a sequence of Laurent polynomials $J_K(n) \in \mathbb{Z}[t,t^{-1}]$ that have no apparent connection to classical knot invariants.  The relationship posited by the $AJ$-conjecture allows us to extract information about one invariant from the other.  For example, if the $AJ$-conjecture is true, then the fact that the $A$-polynomial recognizes the unknot implies the colored Jones function does as well. 

The $AJ$-conjecture has been verified for the trefoil and figure eight knots \cite{Ga}, all torus knots \cite{Hi, Tr3}, some classes of two-bridge knots and pretzel knots \cite{Le, LT}, the knot $7_4$ \cite{GK}, and most cable knots over torus knots, the figure eight knot, and twist knots \cite{Ru, RZ, Tr, Tr1}.  Early methods of proving the $AJ$-conjecture relied on explicit formulas for each invariant.  Recently, methods have been developed in \cite{Tr} to verify the conjecture for cables of knots even when explicit formulas are unknown.  In the current paper, a closed formula for the $A$-polynomial is needed. 

The main result of this paper is an extension of work done by Dennis Ruppe and Anh Tran \cite{Ru, Tr, Tr1}.  \newline

\textit{\textbf{Theorem~\ref{thm:main}} Let $K = \mathfrak{b}(p,m)$ be a two-bridge knot with a reduced alternating diagram $D$ that has $\eta_+$ positive crossings and $\eta_-$ negative crossings.  Let $A_K(M,L)$ be the $A$-polynomial of $K$ and let $C$ be the $(r,2)$-cable of $K$.  Assume the following:}  \newline
\indent \textit{i.  $K$ satisfies the $AJ$-conjecture with an associated degree $p$ annihilator $\alpha_K(t,M,L)$   } \newline
\indent \textit{ii.  $\frac{A_K(M,L)}{L-1}$ is irreducible and $\frac{A_K(M,L)}{L-1} \ne \frac{A_K(M,-L)}{L-1}$.}\newline
\indent \textit{iii.  The matrix $N(-1,M)$ described in Section~\ref{sec:biglemma} has nonzero determinant.} \newline 
\textit{Then for odd integers $r$ with $(r+4\eta_-)(r-4 \eta_+)>0$ the cable knot, $C$, satisfies the~$AJ$~-conjecture. } \newline

The need for condition \textit{iii} will be apparent in Section \ref{sec:biglemma} as it leads to a solution of a system of equations.  \newline  

In Section \ref{sec:prelim}, we review each invariant, state the $AJ$-conjecture, and provide other preliminary material.   In Section~\ref{sec:biglemma}, we establish an important lemma needed to prove the main result in Section~\ref{sec:mainthm}.  Applications of the main result are given in section~\ref{sec:results}.  \newline

The author would like to Anh Tran for his substantial influence on the completion of this paper.  Gratitude is also extended to Cody Armond, Ben Cooper and Steven Landsburg for helpful conversations.

\section{Preliminaries}
\label{sec:prelim}

\subsection{The A-Polynomial}
The $A$-polynomial is the defining polynomial of an algebraic curve in $\mathbb{C}^* \times \mathbb{C}^*$ where $\mathbb{C}^*$ are the nonzero complex numbers \cite{CCGLS}.  Let $K \subset S^3$ be a knot and $\mathcal{M}$ be the complement of a regular neighborhood of $K$.  Then $\mathcal{M}$ is a compact manifold with boundary homeomorphic to a torus, $\partial \mathcal{M} = T$.  Denote by $\textrm{Rep}(T) \subset \prod_{i=1}^n SL_2(\mathbb{C})$ the space of all representations of $\pi_1(T)$ into $SL_2(\mathbb{C})$.  

The fundamental group $\pi_1(T)$ is a free abelian group with two generators.  Let $\{ \mu, \lambda \}$ be the standard basis for $\pi_1(T)$.  Consider the subset $\textrm{Rep}^\Delta(\pi_1(\mathcal{M}))$ of $\textrm{Rep}(\pi_1(\mathcal{M}))$ consisting of upper triangular $ SL_2(\mathbb{C})$ representations.  Set
$$\rho(\mu) =  \begin{pmatrix} m & \star \\ 0  & m^{-1} \end{pmatrix}  \qquad \textrm{and} \qquad \rho(\lambda) = \begin{pmatrix} l & \star \\ 0  & l^{-1} \end{pmatrix} $$
and let $\epsilon : \textrm{Rep}^\Delta(\pi_1(\mathcal{M})) \rightarrow \mathbb{C}^* \times \mathbb{C}^*$ be the eigenvalue map defined by $\epsilon(\rho) = (m,l)$.  Let $Z$ be the Zariski closure of $\epsilon(\textrm{Rep}^\Delta(\pi_1(\mathcal{M})))$ in $\mathbb{C}^* \times \mathbb{C}^*$.  Each of the components of $Z$ are  one dimensional \cite{DG}.  The components are hyper-surfaces and can be cut out by a single polynomial unique up to multiplication by a constant.  The \textbf{$A$-polynomial}, $A_K(m,l)$, is the product of all such defining polynomials.  The $A$-polynomial can be taken to have relatively prime integer coefficients and is well defined up to a unit.  The abelian component of $Z$ will have defining polynomial $l-1$ and thus the $A$-polynomial can be factored as $A_K(m,l) = (l-1)A'_K(m,l)$ \cite{CCGLS}.

\subsection{The Colored Jones Function}

Let $L$ be a link in $S^3$.  The colored Jones function is a quantum link invariant that assigns to each $n \in \mathbb{N}$ a Laurent polynomial $J_L(n) \in \mathbb{Z}[t^{\pm 1}]$.  Like the classic Jones polynomial, $J_L(n)$ can be defined using the Kauffman bracket.  For a complete definition of the colored Jones function, see for example \cite{Le}.

\subsection{The AJ Conjecture}

To set up the $AJ$ conjecture, we follow \cite{Ga} and use notation from \cite{Le}.  Given a discrete function $f:\mathbb{N} \rightarrow \mathbb{C}[t^{\pm1}]$, define operators $M$ and $L$ acting on $f$ by
$$M(f)(n) = t^{2n}f(n) \quad \textrm{and} \quad L(f)(n) = f(n+1)$$
It can be seen that these operators satisfy $LM = t^2 ML$.  Let 
$$\mathcal{T}= \mathbb{C}[t^{\pm1}] \left< M^{\pm1},L^{\pm1} \right> / (LM - t^2ML)        $$
This non-commutative ring, $\mathcal{T}$, is called the quantum torus. \newline

If there exists a $P \in \mathcal{T}$ such that $P(f) = 0$, then $P$ is called a \textit{recurrence relation} for $f$.  The \textbf{recurrence ideal} of a discrete function $f$ is a left ideal $A_f$ of $\mathcal{T}$ consisting of recurrence relations for $f$.
$$A_f = \{ P \in \mathcal{T} | P(f) = 0 \}$$

When $A_f \ne \{0\}$, the discrete function $f$ is said to be \textbf{q-holonomic}.  Also, denote by $A_K$ the recurrence ideal of $J_K(n)$. \newline

The ring $\mathcal{T}$ is not a principal ideal domain, so a non-trivial recurrence ideal is not guaranteed to have a single generator.  Fix this by localizing at Laurent polynomials in $M$ \cite{Ga}.  The process of localization works since $\mathcal{T}$ satisfies the Ore condition.  The resulting ring, $\widetilde{\mathcal{T}}$, has a Euclidean algorithm and is therefore a principal ideal domain.  Let $\mathbb{C}[t^{\pm1}](M)$ be the fraction field of the polynomial ring $\mathbb{C}[t^{\pm1}][M]$.  Then,  
$$\widetilde{\mathcal{T}} = \left\{ \sum_{i \in \mathbb{Z}} a_i(M) L^i \, \middle| \, a_i(M) \in \mathbb{C}[t^{\pm1}](M), \, a_i(M) = 0 \, \textrm{for  almost  every i} \right\}.$$

Extend the recurrence ideal $A_K$ to an ideal $\widetilde{A_K}$ of $\widetilde{\mathcal{T}}$ by $\widetilde{A_f} = \widetilde{\mathcal{T}}A_K$ which can be done since $\mathcal{T}$ embeds as a subring of $\widetilde{\mathcal{T}}$.  The extended ideal will then have a single generator 
$$\widetilde{\alpha_K}(t,M,L) = \sum_{i=0}^n \alpha_i(M)L^i$$

This generator has only positive powers of $M$ and $L$, and the degree of $L$ is assumed to be minimal since it generates the ideal.  The coefficients $\alpha_i(M)$ can be assumed to be co-prime and $\alpha_i(M) \in \mathbb{Z}[t^{\pm1}, M]$.  The operator $\widetilde{\alpha_K}$ is called the \textbf{recurrence polynomial} of $K$ and is defined up to a factor of $\pm t^jM^k$ for $j,k \in \mathbb{Z}$.  Garoufalidis and Le showed that for every knot $K$, $J_K(n)$ satisfies a nontrivial recurrence relation \cite{GL}. \newline

Let $f,g \in \mathbb{C}[M,L]$.  If the quotient $f/g$ does not depend on $L$, then $f$ and $g$ are said to be \textit{$M$-essentially equivalent}, denoted $f \stackrel{M}{=}g$.  In other words, the functions $f$ and $g$ are equal up to a factor only depending on $M$.  Let $\epsilon$ be the map were the substitution $t = -1$ is made.  \newline

\textbf{The AJ Conjecture} \cite{Ga}: If $A_K(M,L)$ is the $A$-polynomial of a knot $K$ and $\widetilde{\alpha_K}(t,M,L)$ is as above, then $\epsilon(\widetilde{\alpha_K}) \stackrel{M}{=} A_K(M,L)$

\subsection{Cabling Formulas and the Resultant}
Let $r,s \in \mathbb{Z}$ with greatest common divisor $d$.  The $(r,s)$-cable, $C$, of a zero framed knot $K$ is the link formed by taking $d$ parallel copies of the $(\frac{r}{d},\frac{s}{d})$ curve on the torus boundary of a tubular neighborhood of $K$.  Homologically, one can think of this curve as $\frac{r}{d}$ times the meridian and $\frac{s}{d}$ times the longitude on the torus boundary.  Notice that if $r$ and $s$ are relatively prime, $C$ is a knot. \newline

Let $A_K(M,L)$ be the $A$-polynomial of a knot $K$ in $S^3$ and let $C$ be the $(r,2)$-cable knot of $K$.  Assume that $A_K'(M,L)$ is irreducible.  Write $A_C(M,L)$ in terms of $A_K(M,L)$ using the following cabling formula given by Ni and Zhang, c.f. \cite{Ru}.  Let 
\begin{displaymath}
   F_r(M,L) = \left\{
     \begin{array}{lr}
       M^{2r}L + 1 & \textrm{if} \quad r > 0\\
       L + M^{-2r} & \textrm{if}  \quad r < 0
     \end{array}
   \right.
\end{displaymath} 
Then,
$$A_C(M,L) = (L-1) F_r(M,L) \textrm{Res}_\lambda \left( A_K'(M^2,\lambda), \lambda^2 - L \right)$$

where Res$_\lambda$ is the resultant defined below.  Use the fact that the product of the resultants is the resultant of the product and $\textrm{Res}_\lambda(\lambda -1,\lambda^2 -1) = L-1$  to rewrite this formula as
\begin{equation} \label{eq:2.1}
A_C(M,L) =  F_r(M,L) \textrm{Res}_\lambda \left( A_K(M^2,\lambda), \lambda^2 - L \right)
\end{equation}

Let $\mathbb{F}$ be a field.  Let $f(x) = \sum_{i=0}^n f_i x^i$ and $g(x) = \sum_{i=0}^m g_i x^i$ be polynomials in $\mathbb{F}[x]$.  The \textbf{resultant} of $f$ and $g$ is the determinant of the following $(m+n) \times (m+n)$ matrix.

$$\textrm{Res}(f,g) = \begin{vmatrix} f_0 & 0 & \hdots & 0 & g_0 & 0 & \hdots & 0 \\
f_1 & f_0 &  & 0 & g_1 & g_0 &  & 0 \\
f_2 & f_1 & \ddots & 0 & g_2 & g_1 & \ddots & 0 \\
\vdots & \vdots & \ddots & \vdots & \vdots & \vdots & \ddots & \vdots \\
f_n & f_{n-1} & \ddots & f_0 & g_k & \vdots & \ddots & g_0 \\
0 & f_n & \ddots & f_1 & g_{k+1} & g_k &\ddots & g_1 \\
\vdots & \vdots & \ddots & \vdots & \vdots & \vdots & \ddots & \vdots \\
0 & \hdots & f_n & f_{n-1} & 0 &  \hdots & g_m & g_{m-1} \\
0 & \hdots & 0 & f_n & 0 & \hdots & 0 & g_m
\end{vmatrix}$$

The coefficients of $f(x)$ constitute the first $m$ columns while the coefficients of $g(x)$ are in the last $n$ columns. \newline

There is also a formula relating the colored Jones polynomial of the cable to a subsequence of the colored Jones polynomial of $K$.  Since we only consider the $(r,2)$-cable knot $C$, the cabling formula given in \cite{RZ} simplifies as
\begin{equation} \label{eq:2.2}
(M^r L + t^{-2r} M^{-r}) J_C(n) = J_K(2n+1).
\end{equation}

In the next section, we seek a homogeneous annihilator $\widetilde{\beta}(t,M,L)$ of the odd subsequence $J_K(2n+1)$.  With this in hand, the above cabling formula implies that $\widetilde{\beta}(t,M,L)(M^r L + t^{-2r} M^{-r})$ annihilates $J_C(n)$. \newline

\section{The Annihilator}
\label{sec:biglemma}

\begin{lemma} \label{biglemma} Let $\widetilde{\alpha}_K(t,M,L) = \sum_{i = 0}^{d} P_i(t,M) L^i$ be a minimal degree homogeneous recurrence polynomial for $J_K(n)$ of degree $d \ge 2$ such that $\widetilde{\alpha}_K(-1,M,L) \stackrel{M}{=} A_K(M,L)$.  Also assume $deg_L(\widetilde{\alpha}_K) = deg_l(A_K)$.  If the matrix $N$ defined later has nonzero determinant at $t=-1$, then $J_K(2n+1)$ has the homogeneous recurrence polynomial given by 
$$\widetilde{\beta}(t,M,L) = \sum_{i = 0}^{d} (-1)^{i} Q_i(t,M) L^i$$
with
$$Q_i(t,M) = -\rm{det}(A_{i+1})$$
where the $A_{i+1}$ are matrices to be defined below.
\end{lemma}$~$ \newline
\textbf{Proof}: The process is similar to the case of the figure eight knot \cite{Ru}, except here a homogeneous annihilator is considered.  Since $\widetilde{\alpha}_K$ is a homogeneous annihilator, $\widetilde{\alpha}_K(t,M,L) J_K(n) = 0$.  Since $M$ acts on discrete functions as multiplication by $t^{2n}$, change $M$ to $t^{2n}$ yielding
$$\left( \sum_{i = 0}^{d} P_i(t,t^{2n}) L^i \right) J_K(n) = \sum_{i = 0}^{d} P_i(t,t^{2n}) J_K(n+i) = 0$$
For each $0 \le j \le d$, substitute $2n+j+1$ for $n$ yielding $d+1$ equations of the form
$$\sum_{i = 0}^{d} P_i(t,t^{2(2n+j+1)}) J_K(2n + j + 1 + i) = 0$$
A degree $d$ homogeneous recurrence relation of $J_K(2n+1)$ has the form 
$$ \left( \sum_{i = 0}^{d} Q_i(t,t^{2n}) L^i \right) J_K(2n+1) = \sum_{i = 0}^{d} Q_i(t,t^{2n})   J_K(2(n+i)+1) = 0$$
To construct $\widetilde{\beta} $ solve 
\begin{equation} \label{eq:3.1}
\sum_{j=0}^{d} c_j \sum_{i = 0}^{d} P_i(t,t^{4n+2j+2}) J_K(2n + i+j+1) - \sum_{i = 0}^{d} Q_i(t,t^{2n})   J_K(2n+2i+1)=0
\end{equation}

Setting the coefficients of each $J_K(m)$ equal to zero, where $m$ satisfies $2n+1 \le m \le 2n + 2d + 1$, gives $2d+1$ equations in the $2d+2$ unknowns $c_0, c_1, ..., c_d, Q_1, Q_2, ... Q_d$.  Form a $(2d+1) \times (2d+2)$ matrix, $C$, where the columns are labeled by the unknowns in the order given above. In the following matrix,  $P_i(r)$ is used as shorthand for $P_i(t,t^{4 n + r})$.

\[
C = 
\left[ \begin{smallmatrix}
P_0(2)&0&0&\cdots&0&\cdots&0&0&-1&0&0&0&\cdots&0&0 \\
P_1(2)&P_0(4)&0&\cdots&0&\cdots&0&0&0&0&0&0&\cdots&0&0 \\
P_2(2)&P_1(4)&P_0(6)&\cdots&0&\cdots&0&0&0&-1&0&0&\cdots&0&0 \\
P_3(2)&P_2(4)&P_1(6)&\ddots&0&\cdots&0&0&0&0&0&0&\cdots&0&0 \\
P_4(2)&P_3(4)&P_2(6)&\ddots&0&\cdots&0&0&0&0&-1&0&\cdots&0&0 \\
P_5(2)&P_4(4)&P_3(6)&\ddots&0&\cdots&0&0&0&0&0&0&\cdots&0&0 \\
\vdots&\vdots&\vdots& \ddots&\vdots& & \vdots&\vdots& \vdots&\vdots&\vdots&\vdots&\vdots&\vdots&\vdots \\
P_{d-1}(2)&P_{d-2}(4)&P_{d-3}(6)&\cdots&0&\cdots&0&0&0&0&0&0&\cdots&0&0 \\
P_d(2)&P_{d-1}(4)&P_{d-2}(6)&\cdots&P_0(2k+2)& &0&0&0&0&0&0&\cdots&0&0\\
0&P_d(4)&P_{d-1}(6)&\cdots&P_1(2k+2)&\ddots&0&0&0&0&0&0&\cdots&0&0 \\
0&0&P_d(6)&\cdots&P_2(2k+2)&\ddots&0&0&0&0&0&0&\cdots&0&0 \\
0&0&0&\cdots&P_3(2k+2)&\ddots&0&0&0&0&0&0&\cdots&0&0 \\
\vdots&\vdots&\vdots&\cdots&\vdots& & \vdots&\vdots&\vdots&\vdots&\vdots&\vdots&\vdots&\vdots&\vdots \\
0&0&0&\cdots&P_{d-1}(2k+2)&\ddots&P_k(2d)&P_{k-2}(2d+2)&0&0&0&0&\cdots&0&0 \\
0&0&0&\cdots&P_d(2k+2)&\ddots&P_{k+1}(2d)&P_k(2d+2)&0&0&0&0&\cdots&0&0 \\
0&0&0& \cdots & 0  & \ddots&P_{k+2}(2d)&P_{k+1}(2d+2)&0&0&0&0&\cdots&0&0\\
\vdots&\vdots&\vdots&& \vdots& &\vdots&\vdots&\vdots&\vdots&\vdots&\vdots& &-1&0 \\
0&0&0&\cdots&0&\cdots &P_{d}(2d)&P_{d-1}(2d+2)&0&0&0&0&\cdots&0&0\\
0&0&0&\cdots&0&\cdots& 0&P_d(2d+2)&0&0&0&0&\cdots&0&-1
\end{smallmatrix}\right]
\]

If $<c_0, ...,c_d,Q_0,...,Q_d>^T = \vec{X} \in \mathbb{R}^{2d+2}$ is the vector of unknowns, then \eqref{eq:3.1} is equivalent to $C \vec{X} = \vec{0}$.  Because $C$ is not a square matrix, Cramer's rule can not be applied.  However, a slight modification can be made to get a square matrix.  \newline

Let $D$ be the $(2d+1) \times (2d+1) $ matrix obtained from $C$ by removing the last column.  Let $\vec{Y} \in \mathbb{R}^{2d +1}$ be the vector obtained from $\vec{X}$ by removing the last entry $Q_d$.  Finally, let $\vec{b} =<0,0,...,0,Q_d>^T \, \in \mathbb{R}^{2d+1}$.  It can be seen that the equation $C \vec{X} = \vec{0}$ is equivalent to $D \vec{Y} = \vec{b}$.  \newline

Assume $\textrm{det}(D) \ne 0$, then by Cramer's rule 
$$c_{k} = \frac{ \textrm{det}(D_{k+1}) }{ \textrm{det}(D)} \qquad \textrm{and} \qquad Q_k =\frac{ \textrm{det}(D_{d+k+2}) }{ \textrm{det}(D)} $$
where $D_j$ is the same as the matrix $D$ except the $j$-th column is replaced by $\vec{b}$.  At this point it is convenient to make the choice $Q_d(t,M) = \textrm{det}(D)$.  In order to verify the AJ-conjecture later, it will be helpful to simplify these determinants.  Let us first examine $\textrm{det}(D)$.  The last $d$ columns of $D$ have only one nonzero entry(that is, $-1$).  Successive row expansion along each of the last $d$ columns yields, 
\begin{equation} \label{eq:3.2}
\textrm{det}(D)  =  \textrm{det}(X)
\end{equation}
where $X$ is a $(d+1) \times (d+1)$ matrix made up of the first $d+1$ columns of $D$ with the odd rows removed except the $2d + 1$ row.  There is no factor of $-1$ in \eqref{eq:3.2} since the column expansion always results in an even number of $-1$'s.  Now simplify a little further to get  
$$ \textrm{det}(X) =  P_d(t,t^{4n+2d}) P_d(t,t^{4n + 2d + 2}) \textrm{det}(N)$$
where $N$ is the upper left $(d-1) \times (d-1)$ block sub-matrix of $X$.  It is the matrix $N$ that will be important in verifying the AJ-conjecture.  \newline

\textbf{Example}: Let $d = 3$.  Then we have

$$ D = \begin{pmatrix} P_0(2) & 0 & 0 & 0 & -1 & 0 & 0  \\
P_1(2) & P_0(4) & 0 & 0 & 0 & 0&0 \\
P_2(2) & P_1(4) & P_0(6) & 0 & 0 & -1 & 0  \\
P_3(2) & P_2(4) & P_1(6) & P_0(8) & 0 & 0 & 0 \\
0 & P_3(4) & P_2(6) & P_1(8) &0&0&-1\\
0&0&P_3(6) & P_2(8) & 0&0&0 \\
0&0&0&P_3(8)&0&0&0
\end{pmatrix} $$

$$N = \begin{pmatrix} 
P_1(2) & P_0(4) \\
P_3(2) & P_2(4)
\end{pmatrix} $$

Now simplify the determinants corresponding to the $c_j$ for $0 \le j \le d$.  For example, $c_0 = \frac{\textrm{det}(D_1)}{\textrm{det}(D)}$.  Simplify $\textrm{det}(D_1)$ by expanding along the last $d$ columns to get a factor of $1$ and then expand along the first column.  Recall that the first column of $D_1$ is the vector $\vec{b}$ with only one nonzero entry, $Q_d$, in the last row.  Therefore, 
$$\textrm{det}(D_1) =(-1)^dQ_d(t,M) \cdot \textrm{det}(B_1)$$
where $B_1$ is the $(d+1,1)$-cofactor of $X$.  In general,
$$\textrm{det}(D_j) =  (-1)^{d+j-1}Q_d(t,M) \cdot \textrm{det}(B_j)$$
where $B_j$ is the $(d+1,j)$ cofactor of $X$ and $1 \le j \le d+1$.  This gives 
$$c_j = (-1)^{d+j} \textrm{det}(B_{j+1}) \quad \textrm{for} \quad 0 \le j \le d.$$

Finally, simplify the determinants corresponding to the $Q_j$ for $0 \le j \le d-1$.  For example, $\textrm{det}(D_{d+2})$ simplifies by expanding along the last $d$ columns, the first of which is the vector $\vec{b}$.  This gives 
$$\textrm{det}(D_{d+2}) = - Q_d \cdot \textrm{det}(A_1)$$ 
where $A_1$ is a $(d+1) \times (d+1)$ sub-matrix of $D$.  The cofactor expansion implies that $A_1$ is formed from the first $d+1$ columns of $D$ while removing all the odd rows except the first.   Similarly, 
$$\textrm{det}(D_{d+j+2}) =- Q_d \cdot \textrm{det}(A_j)$$
where $A_j$ is a $(d+1) \times (d+1)$ sub-matrix of $D$.  Form $A_j$ from the first $d+1$ columns of $D$ while removing all the odd rows except the $2j -1$ row. \newline

\textbf{Example}: In the case $d=3$,

$$A_1 = \begin{pmatrix} 
P_0(2) & 0 & 0 & 0 \\
P_1(2) & P_0(4) & 0 & 0 \\
P_3(2) & P_2(4) & P_1(6) & P_0(8) \\
0 & 0 & P_3(6) & P_2(8)
\end{pmatrix} \qquad  B_1 = \begin{pmatrix} P_0(4) & 0 & 0 \\
P_2(4) & P_1(6) & P_0(8) \\
0 & P_3(6) & P_2(8) \end{pmatrix}$$

$$ A_2 = \begin{pmatrix} 
P_1(2) & P_0(4) & 0 & 0 \\
P_2(2) & P_1(4) & P_0(6) & 0 \\
P_3(2) & P_2(4) & P_1(6) & P_0(8) \\
0 & 0 & P_3(6) & P_2(8)
\end{pmatrix} \qquad B_2 = \begin{pmatrix} P_1(2) & 0 & 0 \\
P_3(2) & P_1(6) & P_0(8)\\
0 & P_3(6) & P_2(8) \end{pmatrix} $$
  
$$A_3 = \begin{pmatrix} 
P_1(2) & P_0(4) & 0 & 0 \\
P_3(2) & P_2(4) & P_1(6) & P_0(8) \\
0 & P_3(4) & P_2(6) & P_1(8) \\
0 & 0 & P_3(6) & P_2(8)
\end{pmatrix} \qquad B_3 = \begin{pmatrix} P_1(2) & P_0(4) & 0\\
P_3(2) & P_2(4) & P_0(8) \\
0 & 0 & P_3(8) \end{pmatrix}$$ $~$
\newline
\newline

To complete the proof of the lemma, it remains to show that
$$\widetilde{\beta}(t,M,L) = \sum_{i = 0}^{d} (-1)^{i} Q_i(t,M) L^i$$
defines a nontrivial operator which can be accomplished with the following claim. \newline

\textbf{Claim}: $\textrm{det}(B_{j+1}(-1,M)) =  P_j(-1,M^2) \cdot \textrm{det}(N(-1,M))$ where $N$ is the $(d-1) \times (d-1) $ matrix defined above. \newline

\textbf{Proof of Claim}:  Break the claim into two cases; $ d$ is even or $d$ is odd.  The case where $d$ is odd is shown, the even case is analogous.  All calculations are done at $t = -1$ and $ P_k $ is shorthand for $P_k(-1,M^2)$. \newline

When $d$ is odd, the matrix $N$ is formed from the two vectors
$$\vec{v} = \begin{pmatrix} P_1 \\ P_3 \\ \vdots \\P_{d-2} \\P_d \\ 0 \\ \vdots \\ 0\end{pmatrix} \quad  \textrm{and} \quad \vec{w} = \begin{pmatrix} P_0 \\ P_2 \\ \vdots \\P_{d-3} \\ P_{d-1} \\ 0 \\ \vdots \\ 0 \end{pmatrix}$$
where $\vec{v}, \vec{w} \in \mathbb{R}^{d-1}$.  For example, when $d=7$

$$\vec{v} = \begin{pmatrix} P_1 \\ P_3 \\ P_5 \\ P_7 \\ 0 \\ 0 \end{pmatrix} \quad  \textrm{and} \quad \vec{w} = \begin{pmatrix} P_0 \\ P_2 \\ P_4 \\ P_6 \\ 0 \\ 0 \end{pmatrix}$$

Let $S$ be the $(d-1) \times (d-1)$ permutation matrix that shifts each vector entry by the permutation $(1 \, 2 \, 3 \, 4 \dots d-1)$.  Then in vector notation, $N$ is the $(d-1) \times (d-1)$ matrix given by 
$$N = \begin{pmatrix}  \vec{v} & \vec{w} & \vec{Sv} & \vec{Sw} & \hdots & \vec{S^{\frac{d-3}{2}}v} & \vec{S^{\frac{d-3}{2}}w} \end{pmatrix}.$$

Let $i: \mathbb{R}^{d-1} \rightarrow \mathbb{R}^d$ be inclusion.  Then $i(\vec{v}) $ is almost the same as $\vec{v}$ except there is an extra zero in the last row of $i(\vec{v})$.  Let $T$ be the permutation matrix that shifts each vector entry by the permutation $(1 \, 2 \, 3 \, 4 \dots d)$.  Form any of the matrices $B_j$ by choosing $d$ vectors (all but the $j$-th vector) from the following set of $d+1$ vectors; 
$$\{ \vec{i(v)} , \vec{i(w)}, \vec{Ti(v)}, \vec{Ti(w)}, \cdots, \vec{T^{\frac{d-1}{2}}i(v)}, \vec{T^{\frac{d-1}{2}}i(w)} \}$$ 

For example, 
$$B_{d+1} =  \begin{pmatrix}  \vec{i(v)} & \vec{i(w)} & \vec{Ti(v)} & \vec{Ti(w)} & \hdots & \vec{T^{\frac{d-3}{2}}i(v)} & \vec{T^{\frac{d-3}{2}}i(w)} &\vec{T^{\frac{d-1}{2}}i(v)} \end{pmatrix}.$$

Notice that $\textrm{det}(B_{d+1}) =  P_d \cdot \textrm{det}(N)$.  This is clear since $B_{d+1}$ is a block matrix.  To show the claim for the rest of the $B_j$, consider the equation
$$B_{d+1} \cdot \vec{X} = \vec{T^{\frac{d-1}{2}}i(w)}$$

This equation has the solution
$$\vec{X} = \frac{1}{P_d} \begin{pmatrix} P_0 \\ - P_1 \\ \vdots \\  -P_{d-2} \\ P_{d-1}  \end{pmatrix}$$

The fact that $\vec{X}$ is the solution to this equation follows from the successive shifting of the vectors $\vec{i(v)}$ and $\vec{i(w)}$.  Notice that because $\widetilde{\alpha_K}$ was a minimal degree annihilator,  $P_0 , P_d \ne 0$, implying $\vec{X}$ is defined and nonzero. \newline

By Cramer's rule, 
$$\frac{ P_i}{P_d} = \frac{\textrm{det}(B_{d+1})_{i+1}}{\textrm{det}(B_{d+1})}$$
where $(B_{d+1})_{i+1}$ is the matrix $B_{d+1}$ with the $i+1$ column replaced by $T^{\frac{d-1}{2}}\vec{i(w)}$.  Note that placing the vector $T^{\frac{d-1}{2}}\vec{i(w)}$ in the $i+1$ column and then exchanging appropriate columns cancels the negative signs from the solution $\vec{X}$.  After this replacement and shifting of columns, the matrix $B_{i+1}$ is obtained.  Therefore,

$$\frac{ P_i}{P_d} = \frac{ \textrm{det}(B_{i+1})}{\textrm{det}(B_{d+1}) }  =  \frac{ \textrm{det}(B_{i+1})}{P_d \cdot \textrm{det}(N) } $$
Conclude that $\textrm{det}(B_{i+1}) = P_i \cdot \textrm{det}(N)$ for $0 \ \le i \le d$ which proves the claim. $\square$ \newline

For general $d$, an important observation is 
$$\textrm{det}(A_i) = \sum_{j+k=2i -1} (-1)^{k+i-1} P_k(2j) \cdot \textrm{det}(B_j)$$ 
where $0 \le k \le d$ and $1 \le j \le d+1$.  This can be calculated directly by expanding along the extra row included in $A_i$ that $B_i$ does not contain.  \newline

The proof of Lemma \ref{biglemma} can be completed by combining the above observation and the previous claim.

$$-Q_i(-1,M) = \textrm{det}(A_{i+1}) = \sum_{j+k=2i +1} (-1)^{k+i} P_k \cdot P_{j-1} \cdot \textrm{det}(N)$$
where $0 \le k \le d$ and $ 1 \le j \le d+1$.  Now shift the $j$ indexing by 1 to get
$$\textrm{det}(A_{i+1}) = \sum_{j+k=2i } (-1)^{k+i} P_k \cdot P_{j} \cdot \textrm{det}(N)$$
where $0 \le k \le d$ and $ 0 \le j \le d$. But then,
\begin{equation} \label{eq:3.3}
(-1)^{i+1} Q_i(-1,M) = (-1)^{2i} \sum_{j+k=2i } (-1)^{k} P_k \cdot P_{j} \cdot \textrm{det}(N).
\end{equation}

The fact that $\widetilde{\alpha}_K(t,M,L) = \sum_{i = 0}^{d} P_i(t,M) L^i$ is a minimal degree annihilator and $deg_L(\widetilde{\alpha}_K) = deg_l(A_K)$ implies that $P_0(-1,M)$ and $P_d(-1,M)$ are nonzero.  Since the condition $\textrm{det}(N(-1,M))\ne0$ has been assumed,  

$$-Q_0(-1,M) = \sum_{j+k=0} (-1)^k P_k(-1,M^2) \cdot P_j(-1,M^2) \cdot \textrm{det}(N(-1,M))$$
$$ = P_0(-1,M^2) \cdot P_0(-1,M^2)  \cdot \textrm{det}(N(-1,M)) \ne 0$$
and
$$-Q_d(-1,M) = \sum_{j+k=2d} (-1)^k P_k(-1,M^2) \cdot P_j(-1,M^2) \cdot \textrm{det}(N(-1,M))$$
$$ = P_d(-1,M^2) \cdot P_d(-1,M^2)  \cdot \textrm{det}(N(-1,M)) \ne 0.$$

Therefore, the operator $\widetilde{\beta}(t,M,L) = \sum_{i = 0}^{d} (-1)^{i} Q_i(t,M) L^i $ is a nontrivial degree $d$ annihilator of the odd sequence $J_K(2n+1)$ completing the proof of Lemma \ref{biglemma}.   $\square$  \newline

Until now, is has been assumed that $\textrm{det}(N) \ne 0$ when evaluated at $t = -1$.  It would be nice if this determinant was nonzero in general.  Generically though, there are cases where it is clearly zero.  For example, assume that a knot $K$ has associated annihilator 
$$\widetilde{\alpha_K}(t,M,L) = \sum_{i=0}^{r} P_{2i}(t,M) \cdot L^{2i}$$
for some positive integer $r$.  This annihilator has only even power of $L$ and thus each column of $N$ corresponding to the odd powers of $L$ is the zero column.  Then clearly $\textrm{det}(N) = 0$.  Therefore, knots where $A_K(M,L)$ has at least one nonzero odd degree $L$ coefficient should be considered to apply Lemma \ref{biglemma}.  \newline

\textbf{Remark}:  For twist knots $K_m$, Anh Tran uses skein theory to show that a variant of $N$ has a nonzero determinant, \cite{Tr}.  In fact, he shows that under a suitable basis, the columns of $N$ are linearly independent.  It should be noted that other authors have encountered this matrix obstruction while proving the $AJ$-conjecture for cables of knots.

\section{Verifying the AJ conjecture}
\label{sec:mainthm}

Let $K$ be a knot that satisfies the $AJ$-conjecture with homogeneous annihilator $\widetilde{\alpha_K}(t,M,L)$ such that $deg_L(\widetilde{\alpha}_K) = deg_l(A_K)$.  Let $C$ be the $(r,2)$ cable knot of $K$.  In the previous section, we found the annihilator $\widetilde{\beta}$ of the colored Jones function $J_K(2n+1)$.  From the cabling formula \eqref{eq:2.2} for the colored Jones function, we have an annihilator of $J_C(n)$ given by

$$\widetilde{\beta}(t,M,L)(M^r L + t^{-2r} M^{-r})$$

Recall \eqref{eq:2.1} the $A$-polynomial of the $(r,2)$-cable is given in terms of the resultant.
$$A_C(M,L) = F_r(M,L) \textrm{Res}_\lambda \left(A_K(M^2,\lambda), \lambda^2 -L\right)$$
where $F_r(M,L) = M^{2r}L +1$ if $r > 0$ and $L + M^{-2r}$ if $r < 0$.  \newline

Therefore, to verify the $AJ$-conjecture it must be shown that
$$\widetilde{\beta}(-1,M,L) = R(M) \textrm{Res}_\lambda \left(A_K(M^2,\lambda), \lambda^2 -L\right)$$
where $R(M)$ is some function of $M$.  \newline

\begin{lemma} \label{lemma:res} Let $P(L,M) = \sum_{i=0}^d P_i(M)L^i$ be a degree $d$ polynomial.  Then,
$$\textrm{Res}_\lambda(P(M^2,\lambda),\lambda^2 - L) =  \sum_{i=0}^d \left( \sum_{k+j =2i } (-1)^k P_k(M^2) P_j(M^2) \right) L^i $$
\end{lemma}
\noindent \textbf{Proof}:  The proof uses induction on the $L$ degree of $P(L,M)$ and the fact that 
$$\textrm{Res}_\lambda(P(M^2,\lambda),\lambda^2 - L) = P(M^2, \sqrt{L}) \cdot P(M^2, -\sqrt{L}).$$  

When $d=2$, calculate the determinant of the resultant matrix to get 
$$\textrm{Res}_\lambda(P(M^2,\lambda),\lambda^2 - L) = {P_2}^2(M^2)L^2 + (2P_0(M^2)P_2(M^2) - {P_1}^2(M^2))L + {P_0}^2(M^2).$$

Now assume that the statement holds for degree less than $d$ and let $P(M,L)$ be a degree $d$ polynomial in $L$.  Then 
$$\textrm{Res}_\lambda(P(M^2,\lambda),\lambda^2 - L) = P(M^2, \sqrt{L}) \cdot P(M^2, -\sqrt{L})$$
$$=(-1)^d{P_d}^2L^d + P_d \sqrt{L^d} \left( \sum_{i=0}^{d-1}(-1)^d  P_i \sqrt{L^i} + P_i (-\sqrt{L})^i \right) + \sum_{i=0}^{d-1} P_i \sqrt{L}^i \cdot \sum_{i=0}^{d-1} P_i \sqrt{-L}^i$$
$$=(-1)^d{P_d}^2 L^d + 2(-1)^d P_d \sqrt{L^d} \left( \sum_{\substack{ i=0 \\i \, \equiv \, d \,\textrm{mod} \, 2}}^{d-1} P_i \sqrt{L^i} \right) + \textrm{Res}_\lambda(\sum_{i=0}^{d-1} P_i \lambda^i , \lambda^2 - L)$$

$$=(-1)^d{P_d}^2 L^d + 2 \left( \sum_{\substack{ i=0 \\i \, \equiv \, d \,\textrm{mod} \, 2}}^{d-1} (-1)^d P_d P_i \sqrt{L^{i+d}} \right) +  \sum_{i=0}^{d-1} \left( \sum_{k+j =2i } (-1)^k P_k P_j \right) L^i$$

$$=  \sum_{i=0}^d \left( \sum_{k+j =2i } (-1)^k P_k P_j \right) L^i  \qed $$

Now use the set up from the previous section to finish verifying the $AJ$-conjecture.  Recall \eqref{eq:3.3}, 
$$(-1)^{i+1} Q_i(-1,M) = (-1)^{2i} \sum_{j+k=2i } (-1)^{k} P_k \cdot P_{j} \cdot \textrm{det}(N).$$
Therefore, by Lemma \ref{lemma:res}, 
\begin{eqnarray*}
\widetilde{\beta}(-1,M,L) &=& \sum_{i = 0}^{d} (-1)^{i} Q_i(t,M) L^i\\
&=& -\sum_{i = 0}^{d} (-1)^{2i} \sum_{j+k=2i } (-1)^{k} P_k \cdot P_{j} \cdot \textrm{det}(N) \cdot L^i\\
&=&- \textrm{det}(N) \sum_{i = 0}^{d}  \sum_{j+k=2i } (-1)^{k} P_k \cdot P_{j} \cdot L^i\\
&=&-\textrm{det}(N) \cdot \textrm{Res}_\lambda(\widetilde{\alpha_K}(-1, M^2,\lambda),\lambda^2 - L)\\
&=&-\textrm{det}(N) \cdot \textrm{Res}_\lambda\left(R(M)\cdot A_K( M^2,\lambda),\lambda^2 - L\right)\\
&=&-\textrm{det}(N) \cdot R^2(M) \cdot \textrm{Res}_\lambda\left( A_K( M^2,\lambda),\lambda^2 - L\right)
\end{eqnarray*}

The $R(M)$ in the last two expressions arises because $K$ satisfies the $AJ$-conjecture.  Therefore, it has been shown that this annihilator of the cable knot $C$ evaluated at $t=-1$ is $M$-essentially equivalent to the  $A$-polynomial $A_C(L,M)$ when $\textrm{det}(N) \ne 0$.  The only thing left to show is that there does not exist an annihilator of lower $L$ degree than $\widetilde{\beta}(t,M,L)(M^r L + t^{-2r} M^{-r})$.  To this end, it is sufficient to put some restrictions on the $A$-polynomial and then use following lemmas and propositions.

\begin{lemma}\label{lemma:square} If $P(M,L) \in \mathbb{C}[M,L]$ is an irreducible polynomial that contains only even powers of $M$, then $P(M^2,L)$ is irreducible over $\mathbb{C}[M,L]$.
\end{lemma}
\textbf{Proof}:  Let $u=M^2$.  Since $P(M,L)$ contains only even powers of $M$, $P(u,L) \in \mathbb{C}[M,L]$.  Assume for contradiction that $P(u^2,L)$ is reducible in $\mathbb{C}[M,L]$.  By symmetry, if $h(M,L)$ is a factor of $P(u^2,L)$ then so is $h(-M,L)$.  Since $P(u,L)$ is irreducible, $P(u^2,L) = h(M,L)h(-M,L)$ where $h(M,L)$ is irreducible in $\mathbb{C}[M,L]$.  If $h(M,L)$ contains a term of odd degree in $M$, let $d$ be the smallest of those degrees.  Then $P(u^2,L)$ contains an term of degree $d$ in $u$, a contradiction.  This means that all the terms in $h(M,L)$ are of even degree in $M$.  Therefore, $P(u,L) = h(u,L)h(-u,L)$ is a polynomial factorization, contradicting the assumption that $P(M,L)$ is irreducible.

\begin{lemma}[Tran \cite{Tr}] \label{lemma:deg} Let $P(M,L) \in \mathbb{C}[M,L]$ be an irreducible polynomial with $P(M,L) \ne P(M,-L)$.  Then, 

$$R_K(M,L) = \textrm{Res}_\lambda(P( M,\lambda),\lambda^2 - L) \in \mathbb{C}[M,L]$$
is irreducible and has $L$ degree equal to that of $P(M,L)$. \newline
\end{lemma}

\begin{lemma} \label{lemma:degree}
Let $K = \mathfrak{b}(p,m)$ be a two-bridge knot.  Assume $K$ satisfies the $AJ$-conjecture.  Also assume $A_K'(m,l)$ is irreducible in $\mathbb{C}[M,L]$ and has $L$ degree $\frac{p-1}{2}$.  Then $\widetilde{\alpha}_K(t,M,L)$ has $L$ degree $\frac{p+1}{2}$.
\end{lemma}

\textbf{Proof}:  By \cite{Le} Proposition 6.1, $\widetilde{\alpha}_K(t,M,L)$ has $L$ degree at most $\frac{p+1}{2}$.  Since $\widetilde{\alpha}_K(-1,M,L) \stackrel{M}{=} (L-1)A_K'(M,L)$ has $L$ degree $\frac{p+1}{2}$, conclude that $\widetilde{\alpha}_K(t,M,L)$ has $L$ degree $\frac{p+1}{2}$. \qed \newline

\begin{prop}[Tran \cite{Tr1}] \label{prop:greater1} Suppose $K$ is a non-trivial alternating knot.  Then the annihilator of the odd sequence $J_K(2n+1)$ has $L$-degree greater than 1.
\end{prop}

\begin{prop}[Tran \cite{Tr1}] \label{prop:div} For any non-trivial annihilator $\widetilde{\delta}(t,M,L)$ of the odd sequence $J_K(2n+1)$, $\epsilon ( \widetilde{\delta})$ is divisible by $L-1$.
\end{prop}

\begin{prop}[Tran \cite{Tr}]\label{prop:cross}  Suppose $K$ is a non-trivial knot with reduced alternating diagram $D$.  Let $\widetilde{\delta}(t,M,L)$ be the minimal degree annihilator of $J_K(2n+1)$ and $\widetilde{\Delta}_C(t,M,L)$ be the minimal degree annihilator of $J_C(n)$.  Then for odd integers $r$ with $(r+4\eta_-)(r-4 \eta_+)>0$, $\widetilde{\Delta}_C(t,M,L) = \widetilde{\delta}(t,M,L)(L+t^{-2r}M^{-2r})$.
\end{prop}

\begin{theorem}\label{thm:main} Let $K = \mathfrak{b}(p,m)$ be a two-bridge knot with a reduced alternating diagram $D$ that has $\eta_+$ positive crossings and $\eta_-$ negative crossings.  Assume $A_K(M,L)$ has $L$ degree $\frac{p+1}{2}$ and  let $C$ be the $(r,2)$-cable of $K$.  Assume the following properties; \newline
 i.  $K$ satisfies the $AJ$-conjecture and $J_K(n)$ has a minimal degree $d$ homogeneous annihilator $\alpha_K(t,M,L)$ with $\widetilde{\alpha}_K(-1,M,L) \stackrel{M}{=} A_K(M,L)$, where $d = \frac{p+1}{2}$. \newline
 ii.  The $A'$ polynomial, $A_K'(M,L)$ of $K$ is irreducible and $A_K'(M,L) \ne A_K'(M,-L)$. \newline
 iii.  The matrix $N(-1,M)$ described in Section~\ref{sec:biglemma} has nonzero determinant. \newline

Then for odd integers $r$ with $(r+4\eta_-)(r-4 \eta_+)>0$ the cable knot $C$ satisfies the $AJ$-conjecture.  
 
\end{theorem}

\textbf{Proof}:  The proof is similar to that of Theorem 1 given by Anh Tran \cite{Tr}.   Conditions $(i)$ and $(iii)$ allow for the application of Lemma \ref{biglemma} to find a non-trivial operator $\widetilde{\beta}(t,M,L)$ such that 
$ \widetilde{\beta}(t,M,L) J_K(2n+1) = 0.$  From the discussion above, $\epsilon(\widetilde{\beta}) \stackrel{M}{=} (L-1)R_K(M^2,L) $ where $R_K(M^2,L) =  \textrm{Res}(A_K'( M^2,\lambda),\lambda^2 - L)$.  Let $\widetilde{\delta}(t,M,L)$ be the minimal degree annihilator of $J_K(2n+1)$.  Then $\epsilon(\widetilde{\delta})$ left divides $\epsilon(\widetilde{\beta})$.  Since $K$ is alternating Propositions \ref{prop:greater1} and \ref{prop:div} imply $\epsilon(\widetilde{\delta})$ has $L$ degree $\ge 2$ and is divisible by $L-1$.  It is well known that the $A$-polynomial has only even powers of $M$.  Therefore, the resultant, $R_K(M^2,L)$, is irreducible over $\mathbb{C}[M,L]$ by condition $(ii)$ and Lemmas \ref{lemma:square} and \ref{lemma:deg}.  Since $\frac{\epsilon(\widetilde{\delta})}{L-1}$ has $L$ degree at least 1 and $R_K(M^2,L)$ is irreducible, conclude that $\frac{\epsilon(\widetilde{\delta})}{L-1} \stackrel{M}{=} R_K(M^2,L)$.  Let $\widetilde{\Delta}_C(t,M,L)$ be the minimal degree annihilator of $J_C(n)$.  All that is left to show is $A_C(M,L) \stackrel{M}{=} \epsilon(\widetilde{\Delta}_C)$.   Since $(r+4\eta_-)(r-4 \eta_+)>0$, Proposition \ref{prop:cross} implies that $\widetilde{\Delta}_C =   \widetilde{\delta} (L + M^{-2r})$.  By these remarks and the cabling formula \eqref{eq:2.1} for $A$-polynomials,
\begin{eqnarray*}
A_C(M,L) &=& (L-1)R_K(M^2,L) (L + M^{-2r})\\
 &\stackrel{M}{=}& \epsilon(\widetilde{\delta}) (L + M^{-2r})\\
&=& \epsilon(\widetilde{\Delta}_C)
\end{eqnarray*}
\vspace{-1cm}
\begin{flushright} $\square$ \end{flushright}

\section{Results} 
\label{sec:results}

To apply Theorem~\ref{thm:main}, first find a two-bridge knot with reduced alternating diagram that satisfies the $AJ$-conjecture.  In \cite{LT}, conditions were given in order for a knot to satisfy the $AJ$ conjecture.  Specifically, all two-bridge knots for which the SL$_2(\mathbb{C})$ character variety has exactly two irreducible components satisfy the $AJ$ conjecture.  We will apply the main theorem to $(r,2)$-cables of some of these two-bridge knots. \newline

\textbf{The knot $6_2 = \mathfrak{b}(11,3)$}

The $A$-polynomial of the knot $6_2$ (without the $L-1$ factor) is irreducible and is given by \newline
$A_{6_2}'(M,L) = -L^5M^{26}+(3M^{18}+2M^{24}+M^{26}-5M^{22}-5M^{20}-2M^{28}+M^{30})L^4+(-3M^{10}+M^{28}-12M^{18}+3M^{14}+8M^{12}+M^{24}-13M^{16}-3M^{26}+5M^{22}+3M^{20})L^3+(8M^{18}+3M^{10}-12M^{12}+M^2-13M^{14}-3M^4+3M^{16}-3M^{20}+5M^{8}+M^6)L^2+(M^4+3M^{12}-2M^2-5M^{10}+1-5M^{8}+2M^6)L-M^4$ \newline

The degree of $A_{6_2}'(M,L)$ is $5 = \frac{11-1}{2}$.  Lemma \ref{lemma:degree} implies the annihilator is degree 6.  The only condition to verify in order to apply Theorem 4.1 is $\textrm{det}(N) \ne 0$.  Since the degree of $A_{6_2}(M,L)$ is 6, $N$ will be a $5 \times 5$ matrix.  Below, the shorthand $P_i = P_i(M^2)$ is used.  
$$N = \begin{pmatrix} P_1 & P_0 & 0 & 0 &0 \\
P_3 & P_2 & P_1 & P_0&0 \\
P_5 & P_4 & P_3 & P_2&P_1 \\
0 & P_6 & P_5 & P_4&P_3 \\
0&0&0&P_6&P_5\end{pmatrix} $$
\newline

With the aid of Maple the determinant of $N$ is found using the coefficients of $A_{6_2}(M,L)$.  The determinant of $N(-1,M)$ will be a nonzero multiple of the below polynomial and therefore still nonzero. \newline

\includegraphics[scale=0.7]{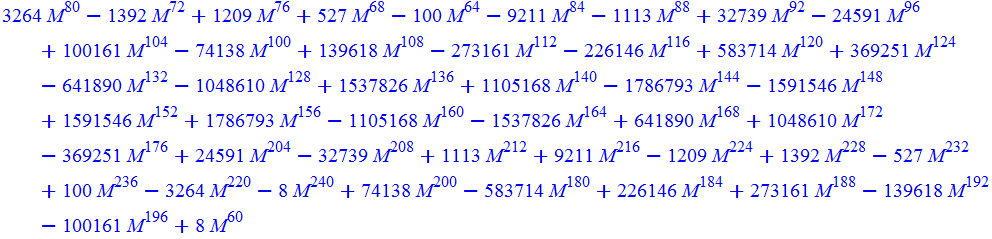} $~$ \newline

The knot $6_2$ has a reduced alternating diagram with four positive and two negative crossings.  Therefore, when $(r-16)(r+ 8) > 0 $, the $(r,2)$-cable knot of $6_2$ will satisfy the $AJ$ conjecture.\newline

Follow the same process to show that the $(r,2)$-cable knots of some two-bridge knots satisfies the $AJ$ conjecture when $r$ is in the proper range.  It would be too cumbersome and take up to much space to list all of the determinants.  However, it has been verified with Maple that each relevant determinant is non-zero.   Below is a table of some two-bridge knots whose cables satisfy the $AJ$ conjecture by Theorem~\ref{thm:main}.  The needed $A$-polynomials were gathered from KnotInfo.

\begin{center}
\begin{tabular}{ c  c}
 \begin{tabular}{ c | c  }
    \hline
    Rolfsen Notation & Two-Bridge Notation \\ \hline
    $6_2$ & $\mathfrak{b}(11,3)$  \\ \hline
   $6_3$ & $\mathfrak{b}(13,5)$  \\ \hline
   $7_3$ & $\mathfrak{b}(13,3)$ \\ \hline
$7_5$ & $\mathfrak{b}(17,7)$  \\ \hline
$7_6$ & $\mathfrak{b}(19,7)$  \\ \hline
$8_2$ & $\mathfrak{b}(17,3)$  \\ \hline
$8_3$ & $\mathfrak{b}(17,13)$  \\ \hline
$8_4$ & $\mathfrak{b}(19,5)$  \\ \hline
$8_6$ & $\mathfrak{b}(23,7)$  \\ \hline
$8_7$ & $\mathfrak{b}(23,5)$  \\ \hline
$8_8$ & $\mathfrak{b}(25,9)$  \\ \hline
$8_9$ & $\mathfrak{b}(25,7)$  \\ \hline
$8_{12}$ & $\mathfrak{b}(29,12)$  \\ \hline
$8_{13}$ & $\mathfrak{b}(29,11)$  \\ \hline
$8_{14}$ & $\mathfrak{b}(31,12)$  \\ \hline
$9_{3}$ & $\mathfrak{b}(19,13)$  \\ \hline

  \end{tabular} &

  \begin{tabular}{ c | c  }
    \hline
    Rolfsen Notation & Two-Bridge Noation \\ \hline
 $9_{4}$ & $\mathfrak{b}(19,13)$  \\ \hline
$9_{5}$ & $\mathfrak{b}(21,5)$  \\ \hline
$9_{7}$ & $\mathfrak{b}(29,13)$  \\ \hline
$9_{8}$ & $\mathfrak{b}(31,11)$  \\ \hline
$9_{9}$ & $\mathfrak{b}(31,9)$  \\ \hline
$9_{13}$ & $\mathfrak{b}(37,27)$  \\ \hline
$9_{14}$ & $\mathfrak{b}(37,14)$  \\ \hline
$9_{15}$ & $\mathfrak{b}(39,16)$  \\ \hline
$9_{17}$ & $\mathfrak{b}(39,14)$  \\ \hline
$9_{18}$ & $\mathfrak{b}(41,17)$  \\ \hline
$9_{19}$ & $\mathfrak{b}(41,16)$  \\ \hline
$9_{20}$ & $\mathfrak{b}(41,15)$  \\ \hline
$9_{21}$ & $\mathfrak{b}(43,18)$  \\ \hline
$9_{26}$ & $\mathfrak{b}(47,18)$  \\ \hline
$9_{27}$ & $\mathfrak{b}(49,19)$  \\ \hline
\\\hline
  \end{tabular}

\end{tabular}
\end{center}

\end{document}